\documentclass[11pt]{article}

\usepackage[a4paper,margin=0.9in]{geometry}
\usepackage{amsmath,amssymb,amsfonts,bm}
\usepackage{booktabs,array,multirow}
\usepackage{xcolor}
\usepackage{float}
\usepackage{pgfplots}
\usepgfplotslibrary{groupplots}

\usepackage{booktabs}
\usepackage{graphicx}

\providecommand{\sci}[2]{%
  \ensuremath{#1\times 10^{#2}}%
}

\usepackage[hidelinks]{hyperref}
\usepackage{tikz}
\usetikzlibrary{calc}
\pgfplotsset{compat=1.18}

\usepackage[
  backend=biber,
  style=numeric,
  sorting=nyt
]{biblatex}

\pgfplotsset{compat=1.18}
\renewcommand{\arraystretch}{0.96}
\makeatletter
\renewcommand\section{\@startsection{section}{1}{\z@}%
  {-2.0ex plus -.5ex minus -.2ex}{1.0ex plus .2ex}%
  {\normalfont\large\bfseries}}
\makeatother

\newcommand{\Om}{\Omega}
\newcommand{\Omone}{\Omega_1}
\newcommand{\Omtwo}{\Omega_2}
\newcommand{\Gam}{\Gamma}
\newcommand{\R}{\mathbb{R}}
\definecolor{revisionred}{RGB}{185,0,0}
\definecolor{additionblue}{RGB}{0,70,175}

\usepackage{authblk}

\title{On the Choice of Continuous and Discontinuous
Lagrange Multipliers in Fictitious-Domain Methods}

\author[1]{Najwa Alshehri}
\author[1]{Shumoua F. Alrzqi%
\thanks{Corresponding author:
\href{mailto:Rizqis@rcjy.edu.sa}{Rizqis@rcjy.edu.sa}}}

\affil[1]{General Studies Department, Jubail Industrial College,
Jubail Industrial City, Saudi Arabia}

\date{}

\begin{document}

\maketitle
\vspace{-1.2cm}

\begin{abstract}
This article examines the Fictitious-Domain approach using Distributed Lagrange Multipliers (FD-DLM), a method that handles elliptic interface problems with separate meshes without requiring the interface to match the background mesh. It compares two stable, low-order FD-DLM discretizations when using both continuous and discontinuous multiplier spaces. The study distinguishes between satisfying the weak coupling constraint per element and the size of the discrepancy between the two discrete primal fields. Numerical experiments indicate that the continuous formulation usually yields smaller local errors and maximum mismatch values, while the discontinuous formulation shows a more regular decrease in mismatch as the mesh is refined and has smaller primal errors.
\end{abstract}

\noindent\textbf{Keywords:} elliptic interface problem; fictitious-domain
method; distributed Lagrange multiplier; unfitted
meshes.

\section{Introduction}

Interface problems arise when distinct materials or physical models meet and interact across an internal boundary, which is commonly seen in composite media, heat transfer, porous transport, and fluid–structure interaction. Discretizing these problems is challenging, particularly when the interface is geometrically complex or dynamic. Under such conditions, a conforming mesh may require frequent remeshing. The Fictitious-Domain Method with Distributed Lagrange Multipliers handles this problem by extending one subdomain problem to a simple background domain. It independently treats the background and immersed domains and enforces agreement through a multiplier distributed over the immersed domain. The main practical value is therefore geometric flexibility without sacrificing a variational treatment of the interface conditions. Continuous-multiplier discretizations for elliptic interface problems have been investigated in \cite{boffi2014mixed}, where mixed fictitious-domain formulations were developed, their well-posedness was established, and good finite element approximation properties were demonstrated. A stable low-order FD-DLM discretization using a discontinuous piecewise-constant multiplier together with a bubble-enriched immersed primal space was examined in ~\cite{alshehri2024unfitted}. A posteriori error estimator for elliptic interface problems was then developed for this formulation in \cite{alshehri2025posteriori}.  More recently, augmented Lagrangian preconditioners designed for FD-DLM elliptic interface problems were developed in \cite{benzi2026augmented}, while multigrid preconditioning for FD-DLM discretizations with both continuous and discontinuous multipliers was studied in \cite{alshehri2025multigrid}. While these studies prove stability and provide effective solvers, they do not determine which multiplier more strongly enforces the discrete coupling condition. This study considers that question. An elementwise-constant discontinuous multiplier evaluates the constraint independently on each immersed element, whereas a continuous multiplier evaluates combinations shared across neighboring elements. This property gives the discontinuous space stronger local weak enforcement. Despite this, it does not determine the magnitude of the actual pointwise or root-mean-square mismatch. We therefore report both types of measures alongside the primal errors, and we distinguish the multiplier degrees of freedom from the size of the complete stable system. 
Section~2 introduces the interface problem and describes the two FD-DLM methods, Section~3 presents the numerical results, and Section~4 provides the paper's conclusions.

\section{Interface problem and FD-DLM formulations}
\label{sec:formulation}

Let $\Om\subset\R^d$, where $d\in\{2,3\}$, denote a bounded domain with a Lipschitz boundary.
The domain is divided into two open subdomains, $\Omone$ and $\Omtwo$, where $\Omega_2$ is completely immersed
in $\Omega$.

\begin{equation*}
    \Om=\Omone\cup\Gam\cup\Omtwo,
\qquad
\partial\Omega=\partial\Omega_1\setminus\Gamma.
\end{equation*}

Suppose the interface $\Gam$ is Lipschitz continuous. The vectors $\bm n_i$ represent the outward unit normal to $\Omega_i$ on $\Gamma$. Let $\beta_i\in L^\infty(\Omega_i)$ be uniformly positive coefficients satisfying $\beta_i\geq\beta_{\min}>0$ almost everywhere, and let $f_i$ be the source terms, then the interface problem is proposed as follows:

\begin{equation}
\begin{aligned}
-\nabla\!\cdot(\beta_i\nabla u_i)&=f_i
&&\text{in }\Om_i,\quad i=1,2,\\
u_1&=u_2
&&\text{on }\Gam,\\
\beta_1\nabla u_1\!\cdot\bm n_1
&=-\beta_2\nabla u_2\!\cdot\bm n_2
&&\text{on }\Gam,\\
u_1&=0
&&\text{on }\partial\Omone\setminus\Gam=\partial\Om .
\end{aligned}
\label{eq:strong}
\end{equation}

For any domain $\omega$, the spaces $L^2(\omega)$ and $H^1(\omega)$ represent the standard
Lebesgue and Sobolev spaces, respectively. The notation $[H^1(\omega)]^*$ refers to the dual space of
$H^1(\omega)$. In the FD-DLM formulation, the functions $u_1$, $f_1$, and $\beta_1$ are extended from $\Omone$ to the entire background domain $\Om$.  We denote the extension of $u_1$ by $u$, which satisfies $u|_{\Omega_2}=u_2$. Similarly, $f$ and $\beta$ are sufficiently smooth extensions of $f_1$ and $\beta_1$. Define

\begin{equation}
V_0=H_0^1(\Om),\quad V_2=H^1(\Omtwo),\quad 
\Lambda=[H^1(\Omtwo)]^* . 
\label{eq:spaces}
\end{equation}

Let $c:\Lambda\times V_2\to\mathbb{R}$ represent the duality pairing,
which is defined by $c(\lambda,v_2):=\langle\lambda,v_2\rangle_{\Lambda,V_2}.$
The continuous FD-DLM formulation aims to find
a solution $(u,u_2,\lambda)\in V_0\times V_2\times\Lambda$ such that

\begin{align}
(\beta\nabla u,\nabla v)_\Om+c(\lambda,v|_{\Omtwo})
  &=(f,v)_\Om &&\forall v\in V_0,\label{eq:cont1}\\
((\beta_2-\beta)\nabla u_2,\nabla v_2)_{\Omtwo}-c(\lambda,v_2)
  &=(f_2-f,v_2)_{\Omtwo} &&\forall v_2\in V_2,\label{eq:cont2}\\
c(\mu,u|_{\Omtwo}-u_2)&=0 &&\forall\mu\in\Lambda .\label{eq:cont3}
\end{align}

Under standard extension and regularity conditions, this variational problem matches the weak formulation of the interface problem~\eqref{eq:strong}. Its well-posedness has been shown in \cite{boffi2014mixed,auricchio2015fictitious}. Let $\mathcal{T}_h$ and $\mathcal{T}_{2,h}$ denote independent, shape-regular quadrilateral meshes of $\Om$ and $\Omtwo$, respectively. For each quadrilateral $K$, $Q_1(K)$ denotes the bilinear
polynomial space of degree at most one in each coordinate direction. $Q_0(K)$ denotes the
space of constant functions. $B(K)$ denotes the biquadratic interior
bubble space whose functions
vanish on $\partial K$. The global $Q_1$ primal spaces are continuous across element boundaries, while $Q_0$ is discontinuous. Let $V_{0,h}:=V_h\cap V_0$. For a stable triple $V_{0,h}\times V_{2,h}\times\Lambda_h$, the discrete problem is aim to find $(u_h,u_{2,h},\lambda_h)\in
V_{0,h}\times V_{2,h}\times\Lambda_h$  such that

\begin{align}
(\beta\nabla u_h,\nabla v_h)_\Om+(\lambda_h,v_h|_{\Omtwo})_{\Omtwo}
 &= (f,v_h)_\Om &&\forall v_h\in V_{0,h},\label{eq:disc1}\\
((\beta_2-\beta)\nabla u_{2,h},\nabla v_{2,h})_{\Omtwo}
 -(\lambda_h,v_{2,h})_{\Omtwo}
 &= (f_2-f,v_{2,h})_{\Omtwo} &&\forall v_{2,h}\in V_{2,h},\label{eq:disc2}\\
(\mu_h,u_h|_{\Omtwo}-u_{2,h})_{\Omtwo}&=0
 &&\forall\mu_h\in\Lambda_h.\label{eq:disc3}
\end{align}

The FD-DLM system is discretized using two stable finite element triples, denoted Element~1 and Element~2.
Element~1 is the continuous stable triple $[Q_1-Q_1-Q_1]$, where $Q_1$  denotes the continuous $Q_1$ multiplier space.
Element~2 is the discontinuous stable triple $[Q_1-(Q_1+B)-Q_0]$, where $Q_0$ consists of all functions
$\mu_h\in L^2(\Omtwo)$ such that
$\mu_h|_K\in Q_0(K)$ for every $K\in\mathcal T_{2,h}$. The stability of Element 1 and Element 2 has been established in \cite{boffi2014mixed,auricchio2015fictitious,alshehri2024unfitted}. It is necessary to include the bubble term in order for Element 2 to be stable. For this reason, the experiment examines fully stable element triples rather than the effect of multiplier continuity separately. Although the discontinuous multiplier reduces the number of unknowns, the need for bubble enrichment increases the total size of the Element 2 system on the refined meshes. The discrete mismatch is given by $e_h:=u_h|_{\Omega_2}-u_{2,h}$. For $Q_0$, selecting the indicator function of any element $K\in\mathcal{T}_{2,h}$ as a test function in equation ~\eqref{eq:disc3}, then

\begin{equation}
\int_K e_h\,dx=0\qquad\forall K\in\mathcal T_{2,h}.
\label{eq:cellmean}
\end{equation}

A cellwise weak identity does not necessarily imply a small pointwise or norm-based discrepancy. This distinction motivates the introduction of complementary
measures of weak-constraint satisfaction and mismatch amplitude:

\begin{equation}
R_K=|K|^{-1}\left|\int_K e_h\,dx\right|,\quad
E_K=\left(|K|^{-1}\int_K e_h^2\,dx\right)^{1/2},\quad
E_K^\infty=\max_{x_q\in\mathcal Q_K}|e_h(x_q)|,
\label{eq:measures}
\end{equation}

here, $\mathcal Q_K$ denotes a common set of quadrature points. The following metrics are reported:
$R_{\max}=\max_KR_K$, $E_{\max}=\max_KE_K$,
$E_{\max}^{\infty}=\max_KE_K^\infty$, and

\begin{equation}
E_{\mathrm{global}}=
\left(|\Omtwo|^{-1}\int_{\Omtwo}e_h^2\,dx\right)^{1/2}.
\label{eq:global}
\end{equation}

Note that $R_{\max}$ quantifies elementwise weak-constraint satisfaction. While, $E_{\max}$, $E_{\max}^{\infty}$, and
$E_{\mathrm{global}}$ measure both the local and global amplitudes of the mismatch between the two discrete primal fields.

\section{Numerical results and interpretation}
\label{sec:results}

Considering $\Om=(-1.4,1.4)^2$, $\Omtwo=\{(x,y):x^2+y^2<1\}$ denotes the immersed unit disk, and $\beta_1=1$.
$\beta_2$ takes values in $\{10, 100\}$, and $f_1 = f_2 = 1$. Let $r^2 = x^2 + y^2$; the exact
solution is given by

\begin{equation}
u_1(r)=\frac{4-r^2}{4},\qquad
u_2(r)=\frac{3\beta_2+1-r^2}{4\beta_2}.
\label{eq:exact}
\end{equation}

The computed values for $h+h_2$ are $0.8879$, $0.4522$, $0.2278$, $0.1142$, and
$0.0572$. The computational meshes are independent, the numerical residuals, evaluated independently, are approximately
$10^{-13}$ to $10^{-12}$.\\

Table~\ref{tab:numerical-summary}(a) reports the primal $L^2$ and $H^1$ errors over the five refinement levels, and Table~\ref{tab:numerical-summary}(b) compares the weak-constraint residual and mismatch measures on the finest mesh. Figure~1 illustrates how the weak-constraint residual and mismatch measures vary under mesh refinement. For both coefficient contrasts and at all five refinement levels, Element~2 always produces smaller $L^2$ errors in both primal fields. For Element~2, $R_{\max}$ varies from about $10^{-10}$ to $10^{-12}$. This matches the cellwise weak identity in equation ~\eqref{eq:cellmean}, which results from algebraic and quadrature errors. On the other hand, Element~1 gives values of the order of $10^{-4}$ to $10^{-3}$. Although there is a clear difference in the weak residual, this does not result in a corresponding difference in the norm-based mismatch measures. On the coarse meshes, Element~1 tends to produce smaller $E_{\max}$ and $E_{\max}^{\infty}$ values, whereas the respective Element~2 quantities decrease more systematically as the mesh is refined. On the finest mesh the local mismatch values obtained from the two formulations differ by less than a factor of 1.2. When $\beta_2=100$, the global mismatch is reduced for Element~2, being $4.060\times10^{-5}$ as compared to $4.296\times10^{-5}$ for Element~1. Over the five refinement levels, the empirical local-mismatch rates lie between $0.63$ and $0.72$ for Element~1 and from $0.92$ to $1.08$ for Element~2. The corresponding global-mismatch rates range from $1.29$ to $1.33$ for Element~1 and from $1.53$ to $1.75$ for Element~2. The empirical rates mentioned above describe the decrease observed in $E_{\max}$, $E_{\max}^{\infty}$, and $E_{\mathrm{global}}$ over the five refinement levels. These rates should not be interpreted as convergence rates for the primal $L^2$ and $H^1$ errors reported in Table~\ref{tab:numerical-summary}(a). The results show that local weak-constraint satisfaction and the agreement of the discrete primal fields are separate properties. The cellwise weak identity controls the element-average mismatch but does not limit its amplitude. This distinction appears since positive and negative contributions may cancel within an immersed element, yet still be detectable by $E_{\max}$ and $E_{\max}^{\infty}$. Because $u|_{\Omega_2}=u_2$, the mismatch is given by $e_h=(u_h-u)|_{\Omega_2}-(u_{2,h}-u_2)$, and as a result the cancellation of the two primal approximation errors may cause $e_h$ to be smaller than each of the individual errors. Because of this, weak residuals, norm-based mismatch measures, and primal errors each capture different parts of discrete coupling, so they should all be considered together. The two approaches also differ in multiplier continuity, how they enrich the immersed primal space, and the overall system size. For this reason, the findings should be seen as features of the full stable triples, not simply as effects of multiplier continuity.

% ============================================================
\begin{table}[H]
\centering
\caption{Numerical comparison of Element 1 and Element 2 stable discretizations}
\label{tab:numerical-summary}

\begingroup
\scriptsize
\renewcommand{\arraystretch}{1.05}
\setlength{\tabcolsep}{3pt}

% ============================================================
% Panel (a): Primal L2 and H1 errors
% ============================================================

\vspace{0.5em}

\textbf{(a) Primal $L^2$ and $H^1$ errors over five refinement levels.}

\vspace{0.3em}

% ------------------------------------------------------------
% L2 errors
% ------------------------------------------------------------

\resizebox{\textwidth}{!}{%
\begin{tabular}{c|c|r|r|r|r|r|r|r|r}
\toprule
\multicolumn{2}{c}{} &
\multicolumn{4}{c}{$\beta_2=10$} &
\multicolumn{4}{c}{$\beta_2=100$}
\\
\midrule

\multicolumn{2}{c|}{} &
\multicolumn{2}{c|}{$\lVert u-u_h\rVert_{L^2(\Omega)}$} &
\multicolumn{2}{c|}{$\lVert u_2-u_{2,h}\rVert_{L^2(\Omega_2)}$} &
\multicolumn{2}{c|}{$\lVert u-u_h\rVert_{L^2(\Omega)}$} &
\multicolumn{2}{c}{$\lVert u_2-u_{2,h}\rVert_{L^2(\Omega_2)}$}
\\
\midrule

Cycle & $h+h_2$
& Element 1 & Element 2
& Element 1 & Element 2
& Element 1 & Element 2
& Element 1 & Element 2
\\
\midrule

0 & 0.8879
& \sci{4.06}{-2} & \sci{2.91}{-2}
& \sci{4.79}{-2} & \sci{3.04}{-2}
& \sci{7.33}{-2} & \sci{4.94}{-2}
& \sci{9.33}{-2} & \sci{6.11}{-2}
\\

1 & 0.4522
& \sci{2.02}{-2} & \sci{1.20}{-2}
& \sci{2.53}{-2} & \sci{1.32}{-2}
& \sci{4.57}{-2} & \sci{3.43}{-2}
& \sci{6.00}{-2} & \sci{4.44}{-2}
\\

2 & 0.2278
& \sci{9.10}{-3} & \sci{2.80}{-3}
& \sci{1.17}{-2} & \sci{3.00}{-3}
& \sci{2.96}{-2} & \sci{2.23}{-2}
& \sci{3.94}{-2} & \sci{2.96}{-2}
\\

3 & 0.1142
& \sci{4.00}{-3} & \sci{1.80}{-3}
& \sci{5.30}{-3} & \sci{2.30}{-3}
& \sci{1.90}{-2} & \sci{1.71}{-2}
& \sci{2.55}{-2} & \sci{2.30}{-2}
\\

4 & 0.0572
& \sci{2.00}{-3} & \sci{1.00}{-3}
& \sci{2.70}{-3} & \sci{1.30}{-3}
& \sci{1.58}{-2} & \sci{1.49}{-2}
& \sci{2.12}{-2} & \sci{2.00}{-2}
\\

\bottomrule
\end{tabular}%
}

\vspace{0.6em}

% ------------------------------------------------------------
% H1 errors
% ------------------------------------------------------------

\resizebox{\textwidth}{!}{%
\begin{tabular}{c|c|r|r|r|r|r|r|r|r}
\toprule
\multicolumn{2}{c}{} &
\multicolumn{4}{c}{$\beta_2=10$} &
\multicolumn{4}{c}{$\beta_2=100$}
\\
\midrule

\multicolumn{2}{c|}{} &
\multicolumn{2}{c|}{$\lVert u-u_h\rVert_{H^1(\Omega)}$} &
\multicolumn{2}{c|}{$\lVert u_2-u_{2,h}\rVert_{H^1(\Omega_2)}$} &
\multicolumn{2}{c|}{$\lVert u-u_h\rVert_{H^1(\Omega)}$} &
\multicolumn{2}{c}{$\lVert u_2-u_{2,h}\rVert_{H^1(\Omega_2)}$}
\\
\midrule

Cycle & $h+h_2$
& Element 1 & Element 2
& Element 1 & Element 2
& Element 1 & Element 2
& Element 1 & Element 2
\\
\midrule

0 & 0.8879
& \sci{2.478}{-1} & \sci{2.362}{-1}
& \sci{5.930}{-2} & \sci{4.690}{-2}
& \sci{2.997}{-1} & \sci{2.604}{-1}
& \sci{1.024}{-1} & \sci{7.390}{-2}
\\

1 & 0.4522
& \sci{1.713}{-1} & \sci{1.634}{-1}
& \sci{3.710}{-2} & \sci{3.310}{-2}
& \sci{2.201}{-1} & \sci{1.934}{-1}
& \sci{7.290}{-2} & \sci{6.100}{-2}
\\

2 & 0.2278
& \sci{1.126}{-1} & \sci{1.041}{-1}
& \sci{2.550}{-2} & \sci{2.280}{-2}
& \sci{1.600}{-1} & \sci{1.302}{-1}
& \sci{5.730}{-2} & \sci{5.150}{-2}
\\

3 & 0.1142
& \sci{7.800}{-2} & \sci{7.070}{-2}
& \sci{1.400}{-2} & \sci{1.550}{-2}
& \sci{1.088}{-1} & \sci{9.690}{-2}
& \sci{4.880}{-2} & \sci{4.770}{-2}
\\

4 & 0.0572
& \sci{5.510}{-2} & \sci{5.040}{-2}
& \sci{9.100}{-3} & \sci{1.150}{-2}
& \sci{8.000}{-2} & \sci{7.430}{-2}
& \sci{4.670}{-2} & \sci{4.620}{-2}
\\

\bottomrule
\end{tabular}%
}

\vspace{0.8em}

% ============================================================
% Panel (b): Finest-mesh measures
% ============================================================

\textbf{(b) Weak-constraint residual and mismatch measures on the
finest mesh, characterized by $h+h_2=0.0572$.}

\vspace{0.3em}

\setlength{\tabcolsep}{5pt}

\begin{tabular}{c|c|r|r|r|r}
\toprule
$\beta_2$
& Method
& $R_{\max}$
& $E_{\max}$
& $E_{\max}^{\infty}$
& $E_{\mathrm{global}}$
\\
\midrule

10
& Element 1
& \sci{4.969}{-4}
& \sci{6.462}{-4}
& \sci{1.564}{-3}
& \sci{4.854}{-5}
\\

10
& Element 2
& \sci{4.012}{-12}
& \sci{7.574}{-4}
& \sci{1.581}{-3}
& \sci{6.669}{-5}
\\

\midrule

100
& Element 1
& \sci{4.552}{-4}
& \sci{5.710}{-4}
& \sci{1.320}{-3}
& \sci{4.296}{-5}
\\

100
& Element 2
& \sci{8.337}{-12}
& \sci{6.389}{-4}
& \sci{1.325}{-3}
& \sci{4.060}{-5}
\\

\bottomrule
\end{tabular}

\endgroup
\end{table}

% ============================================================

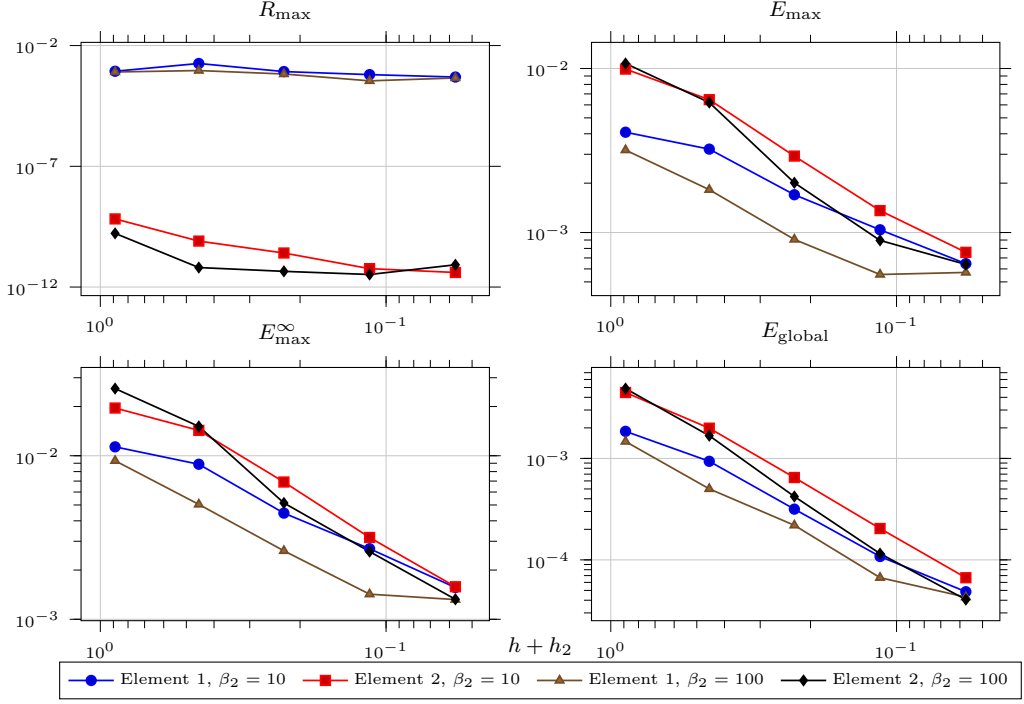
\begin{figure}[H]
\centering
\begingroup
\color{black}

\begin{tikzpicture}
\begin{groupplot}[
    group style={
        group size=2 by 2,
        horizontal sep=1.35cm,
        vertical sep=0.95cm
    },
    width=0.425\textwidth,
    height=0.30\textwidth,
    xmode=log,
    ymode=log,
    x dir=reverse,
    grid=major,
    major grid style={
        draw=gray!45,
        line width=0.25pt
    },
    axis line style={black},
    tick style={black},
    tick align=outside,
    tick label style={
        font=\tiny,
        text=black
    },
    label style={
        font=\scriptsize,
        text=black
    },
    title style={
        font=\scriptsize,
        text=black,
        yshift=-1pt
    },
    every axis plot/.append style={
        line width=0.65pt,
        mark size=1.8pt
    },
    legend style={
        font=\tiny,
        text=black,
        draw=black,
        legend columns=4,
        /tikz/every even column/.append style={
            column sep=4pt
        }
    }
]

%------------------------------------------------
% Weak cell-average residual
%------------------------------------------------
\nextgroupplot[
    title={$R_{\max}$},
    legend to name=commonlegend
]

\addplot+[mark=*]
coordinates {
    (0.8879,8.630585e-4)
    (0.4522,1.811359e-3)
    (0.2278,8.365372e-4)
    (0.1142,6.260754e-4)
    (0.0572,4.968593e-4)
};
\addlegendentry{Element 1, $\beta_2=10$}

\addplot+[mark=square*]
coordinates {
    (0.8879,6.636190e-10)
    (0.4522,8.001199e-11)
    (0.2278,2.568821e-11)
    (0.1142,5.782577e-12)
    (0.0572,4.011700e-12)
};
\addlegendentry{Element 2, $\beta_2=10$}

\addplot+[mark=triangle*]
coordinates {
    (0.8879,8.050487e-4)
    (0.4522,9.232704e-4)
    (0.2278,6.644514e-4)
    (0.1142,3.457890e-4)
    (0.0572,4.551804e-4)
};
\addlegendentry{Element 1, $\beta_2=100$}

\addplot+[mark=diamond*]
coordinates {
    (0.8879,1.672161e-10)
    (0.4522,6.423047e-12)
    (0.2278,4.469823e-12)
    (0.1142,3.310881e-12)
    (0.0572,8.337082e-12)
};
\addlegendentry{Element 2, $\beta_2=100$}

%------------------------------------------------
% Local RMS mismatch
%------------------------------------------------
\nextgroupplot[
    title={$E_{\max}$}
]

\addplot+[mark=*]
coordinates {
    (0.8879,4.088170e-3)
    (0.4522,3.224052e-3)
    (0.2278,1.699998e-3)
    (0.1142,1.039417e-3)
    (0.0572,6.462234e-4)
};

\addplot+[mark=square*]
coordinates {
    (0.8879,9.917498e-3)
    (0.4522,6.457950e-3)
    (0.2278,2.924193e-3)
    (0.1142,1.360870e-3)
    (0.0572,7.574032e-4)
};

\addplot+[mark=triangle*]
coordinates {
    (0.8879,3.182347e-3)
    (0.4522,1.824432e-3)
    (0.2278,9.074243e-4)
    (0.1142,5.547765e-4)
    (0.0572,5.710285e-4)
};

\addplot+[mark=diamond*]
coordinates {
    (0.8879,1.072954e-2)
    (0.4522,6.202897e-3)
    (0.2278,2.010371e-3)
    (0.1142,8.953131e-4)
    (0.0572,6.388717e-4)
};

%------------------------------------------------
% Maximum pointwise mismatch
%------------------------------------------------
\nextgroupplot[
    title={$E_{\max}^{\infty}$}
]

\addplot+[mark=*]
coordinates {
    (0.8879,1.133503e-2)
    (0.4522,8.870870e-3)
    (0.2278,4.458339e-3)
    (0.1142,2.700746e-3)
    (0.0572,1.563751e-3)
};

\addplot+[mark=square*]
coordinates {
    (0.8879,1.957826e-2)
    (0.4522,1.429197e-2)
    (0.2278,6.909254e-3)
    (0.1142,3.163922e-3)
    (0.0572,1.581462e-3)
};

\addplot+[mark=triangle*]
coordinates {
    (0.8879,9.320484e-3)
    (0.4522,5.041145e-3)
    (0.2278,2.621231e-3)
    (0.1142,1.424102e-3)
    (0.0572,1.319600e-3)
};

\addplot+[mark=diamond*]
coordinates {
    (0.8879,2.575170e-2)
    (0.4522,1.511316e-2)
    (0.2278,5.142435e-3)
    (0.1142,2.593363e-3)
    (0.0572,1.324599e-3)
};

%------------------------------------------------
% Global mismatch
%------------------------------------------------
\nextgroupplot[
    title={$E_{\mathrm{global}}$}
]

\addplot+[mark=*]
coordinates {
    (0.8879,1.849268e-3)
    (0.4522,9.367111e-4)
    (0.2278,3.162969e-4)
    (0.1142,1.079197e-4)
    (0.0572,4.853811e-5)
};

\addplot+[mark=square*]
coordinates {
    (0.8879,4.463590e-3)
    (0.4522,1.986496e-3)
    (0.2278,6.483024e-4)
    (0.1142,2.037621e-4)
    (0.0572,6.668747e-5)
};

\addplot+[mark=triangle*]
coordinates {
    (0.8879,1.467210e-3)
    (0.4522,5.010342e-4)
    (0.2278,2.191702e-4)
    (0.1142,6.675060e-5)
    (0.0572,4.296072e-5)
};

\addplot+[mark=diamond*]
coordinates {
    (0.8879,4.880972e-3)
    (0.4522,1.674740e-3)
    (0.2278,4.212593e-4)
    (0.1142,1.152041e-4)
    (0.0572,4.059937e-5)
};

\end{groupplot}

% Common horizontal-axis label
\node[
    font=\scriptsize,
    text=black
] at
($(group c1r2.south east)!0.5!(group c2r2.south west)+(0,-0.34cm)$)
{$h+h_2$};

% Common legend
\node[text=black] at
($(group c1r2.south east)!0.5!(group c2r2.south west)+(0,-0.79cm)$)
{\pgfplotslegendfromname{commonlegend}};

\end{tikzpicture}

\caption{Behavior of the cellwise weak-constraint residual and the
associated local and global mismatch measures over five
mesh-refinement levels.}
\label{fig:convergence}

\endgroup
\end{figure}

\section{Conclusions}
\label{sec:conclusion}

The study considers two stable low-order formulations that use either continuous or discontinuous Lagrange multipliers. Numerical experiments determine that the discontinuous formulation yields near-zero $R_{\max}$ and lower primal $L^2$ errors. On the other hand, the continuous formulation typically produces smaller local mismatches. Therefore, the selection of multipliers should balance the required accuracy with
the size of the system that results.

\printbibliography

\end{document}